\pdfoutput=1
\documentclass[10pt]{article}
\usepackage[a4paper, margin=1in]{geometry}
\usepackage{changepage}
\usepackage{amsmath, amssymb, amsthm}
\usepackage[utf8]{inputenc}
\usepackage[english]{babel}
\usepackage[numbers]{natbib}
\usepackage{url}
\usepackage[usenames]{color}
\usepackage[colorlinks=true]{hyperref}
\usepackage{cleveref}
\usepackage{tabularx}
\usepackage{breqn}
\usepackage[title]{appendix}
\usepackage{xcolor}
\usepackage{comment}
\theoremstyle{plain}

\newtheorem{corollary}{Corollary}[section]
\newtheorem{lemma}{Lemma}[section]
\newtheorem{definition}{Definition}[section]
\newtheorem{theorem}{Theorem}[section]
\newtheorem*{theorem*}{Theorem}
\newtheorem{remark}{Remark}[section]

\crefname{conjecture}{Conjecture}{Conjectures}
\crefname{theorem}{Theorem}{Theorems}
\crefname{theorem*}{Theorem}{Theorems}
\crefname{corollary}{Corollary}{Corollaries}
\crefname{lemma}{Lemma}{Lemmas}
\crefname{proposition}{Proposition}{Propositions}
\crefname{remark}{Remark}{Remarks}
\crefname{note}{Note}{Notes}
\crefname{definition}{Definition}{Definitions}
\crefname{notation}{Notation}{Notations}
\crefname{example}{Example}{Examples}
\crefname{question}{Question}{Questions}
\crefname{section}{\S}{Sections}
\crefname{equation}{Equation}{Equations}
\crefformat{equation}{(eq.~#2#1#3)}

\newcommand{\floor}[1]{\left\lfloor #1 \right\rfloor}

\title{Arithmetic-term representations for the greatest common divisor}
\author{Mihai Prunescu \footnote{Research Center for Logic, Optimization and Security (LOS), Faculty of Mathematics and Computer Science, University of Bucharest, Academiei 14, Bucharest (RO-010014), Romania; Institute of Logic and Data Science, Bucharest, Romania; Simion Stoilow Institute of Mathematics of the Romanian Academy, Research unit 5, P. O. Box 1-764, Bucharest (RO-014700), Romania. E-mail: {\tt mihai.prunescu@imar.ro}, {\tt mihai.prunescu@gmail.com}.},
Joseph M. Shunia \footnote{Wraithwatch, Austin, TX, USA. E-mail: {\tt jshunia@gmail.com}}}
\date{November 2024 \\ \footnotesize{Revised: June 2025}}

\begin{document}

\maketitle

\begin{abstract} \noindent
We construct a new arithmetic-term representation for the function $\gcd(a,b)$. As a byproduct, we also deduce a representation $\gcd(a,b)$ by a modular term in integer arithmetic.
\\[2mm]
\end{abstract}

\section{Introduction}

By $\gcd(a,b)$ we mean the greatest common divisor of the natural numbers $a$ and $b$.

In \cite{mazzanti2002plainbases} Mazzanti proved the following identity:
$$\gcd(a,b) = \floor{\frac{(2^{a^2 b(b+1)} - 2^{a^2 b}) (2^{a^2 b^2} - 1)}{(2^{a^2 b} - 1)(2^{ab^2}-1)2^{a^2 b^2}}} \bmod 2^{ab},$$
for all $a, b \in \mathbb N$ with $a, b \geq 1$. Marchenkov confirmed this formula in \cite{marchenkov2007superposition}, where he presents a slightly different proof.

This identity had been the missing puzzle-piece to show that every Kalmar elementary function is representable as an arithmetic term. Arithmetic terms are defined in \cite{mazzanti2002plainbases} and \cite{marchenkov2007superposition} as terms built up in the language $L = \{+, \dot{-}, \cdot, /, x^y\}$ with variables and constants interpreted in the set of natural numbers. The modified difference $x \dot{-} y$ is defined $0$ if $x < y$. For some applications, see \cite{shunia2023simple, prunescu2024numbertheoreticfunctions}.

As Mazzanti's effective method to compute arithmetic-term representations depends explicitly on the arithmetic-term representation of the greatest common divisor, every simplification of this term has theoretical importance.

Shunia conjectured in \cite{shunia2024elementary} that Mazzanti's formula might be simplified. Here we present a straight-forward construction of a simpler term. Our construction is independent of Mazzanti's, although both of them use Lemma \ref{LemmaBurton}.

\section{Prerequisites}
Under $\mathbb N$ we understand the set of natural numbers including $0$.

We will need the following facts. The first Lemma will be used to make a connection between the value $\gcd(a,b)$ and the number of solutions in natural numbers of a particular linear Diophantine equation.

\begin{lemma}\label{LemmaBurton} (D. Burton, \cite{burtonnumbertheory} Theorem 2-9, page 40)
    The Diophantine linear equation
    $$ax + by = c$$
    has solutions $(x,y) \in \mathbb Z \times \mathbb Z$ if and only if $d = \gcd(a,b) \,|\, c$.
    If $(x_0, y_0) \in \mathbb Z \times \mathbb Z$ is a solution of the equation $ax + by = c$, then all solutions are given by $x = x_0 + (b/d) \cdot t$ and $y = y_0 - (a/d) \cdot t$ where $t \in \mathbb Z$ is arbitrary.
\end{lemma}

The second Lemma is a tool which proved to be useful for representing sequences, whose generating functions are rational, as arithmetic terms.

\begin{lemma}\label{LemmaTermExtraction} (M. Prunescu \& L. Sauras-Altuzarra, \cite{prunescusauras2024representationcrecursive}, Theorem 1) Let $f \in \mathbb N[[z]]$ be a formal series with radius of convergence $R>0$ at $0$,
$$f(z) = \sum _{n \geq 0} s(n) z^n.$$
Let $c, m \in \mathbb N$ be such that $c^{-1} < R$ and for all $n \geq m$, $s(n) < c^{n-2}$. Then for all $n \geq m$,
$$ s(n) = \left \lfloor c^{n^2} f( c^{-n}) \right \rfloor \bmod c^n. $$
\end{lemma}

The third Lemma is an identity in modular arithmetic.

\begin{lemma}\label{LemmaModIdentity}
    Let $A, B, C \in \mathbb Z$ such that $A, B> 0$, $B \neq 0$,  $C \geq 2$, $C \, |\,A$, $B \nmid A$, $B \bmod C = 1$ and $\lfloor A/B \rfloor \bmod C \neq C-1$. Then the following identity holds:
    $$( (-A) \bmod B) \bmod C = 1 + \lfloor A/B \rfloor \bmod C.$$
\end{lemma}

\begin{proof}
    One observes that:
     $$( (-A) \bmod B) \bmod C = ( (-A) - B \cdot \lfloor (-A)/B \rfloor ) \bmod C =$$
     $$ = ( (-A) \bmod C - (B \bmod C) \cdot (\lfloor (-A)/B \rfloor \bmod C ) ) \bmod C =$$
     $$= ( 0 - (-1) \cdot (\lfloor (-A)/B \rfloor \bmod C ) ) \bmod C = ( - (\lfloor (-A)/B \rfloor \bmod C ) ) \bmod C =$$
     $$ = ( - \lfloor (-A) / B \rfloor ) \bmod C. $$
     On the other hand, for $A, B > 0$  such that $B \nmid A$,
     $$- [ (-A) / B ]  = [ A / B ] + 1, $$
     so the quantity to compute is:
     $$ ( \lfloor A / B \rfloor + 1) \bmod C.$$
     Further we recall that for $x \geq 0$ and $C \geq 2$ such that $x \not \equiv ( C-1 ) \mod C$,
     $$(1 + x) \bmod C = 1 + (x \bmod C).$$
     We apply this fact for $x = \lfloor A/B \rfloor $ to get:
     $$1 + \lfloor A / B \rfloor \bmod C.$$
\end{proof}

\section{A family of C-recursive sequences}

\begin{definition}
    C-recursive sequences of order $d$ are sequences $s : \mathbb N \rightarrow \mathbb C$ satisfying a relation of recurrence with constant coefficients:
    $$s(n+d) + a_1 s(n+d-1) + \dots + a_{d-1} s(n+1) + a_d s(n) = 0$$
    for all $n \geq 0$.
\end{definition}

According with Theorem 4.1.1 in \cite{stanley}, and with Theorem 1 in \cite{PetkovsekZakrajsek}, the C-recursive sequences are exactly the sequences $(s(n))$ such that the generating function:
$$f(z) = \sum_{n \geq 0} s(n) z^n$$
is a rational function $A(z)/B(z)$ with $\deg(A) < \deg(B)$.

\begin{definition}
     Let $a, b$ be natural numbers with $a,b \geq 1$ and let $(s_{a,b}(n))$ be the sequence with generating function:
    $$f_{a,b}(z) := \frac{1}{(z^a - 1)(z^b - 1)}.$$
\end{definition}

\begin{lemma}\label{LemmaTermAB}
    For all $a, b \in \mathbb N$, $a, b \geq 1$, the sequence $(s_{a,b}(n))$ consists of natural numbers. In particular,
    $$s_{a,b}(ab) = \gcd(a,b) + 1.$$
\end{lemma}

\begin{proof}
    We observe that:
    $$\frac{1}{(z^a - 1)(z^b - 1)} = (1 + z^a + z^{2a} + z^{3a} + \cdots)(1 + z^b + z^{2b} + z^{3b} + \cdots).$$
    It follows that for all $n \in \mathbb N$, $s_{a,b}(n)$ is indeed a natural number. Moreover, the value of $s_{a,b}(n)$ equals the number of solutions $(x,y) \in \mathbb N \times \mathbb N$ of the equation:
    $$ax + by = n.$$
    For $n = ab$ we observe that $(x_0, y_0) = (b, 0)$ is indeed a solution of the equation:
    $$ax + by = ab.$$
    By Lemma \ref{LemmaBurton}, the set of solutions $(x,y) \in \mathbb Z \times \mathbb Z$ is parametrized by a parameter $t \in \mathbb Z$ and the solutions have the shape:
    \begin{eqnarray*}
        x & = & b - \frac{b}{d} \cdot t, \\
        y & = &  \frac{a}{d} \cdot t,
    \end{eqnarray*}
    where $d = \gcd(a,b)$.
    As only the non-negative solutions count, $t$ must fulfill:
    $$ b \geq \frac{b}{d} \cdot t \,\, \wedge \,\, t \geq 0.$$
    So $0 \leq t \leq d = \gcd(a,b)$, hence there are $\gcd(a,b) + 1$ many solutions in $\mathbb N \times \mathbb N$.
\end{proof}

\begin{lemma}\label{LemmaUniformBound}
If $a, b \geq 1$ are natural numbers and at least one of them is bigger than $1$, then for all $n \in \mathbb N$ one has:
$$s_{a,b} (n) \leq s_{1,1}(n) = n+1.$$
One has $s_{a,b}(n) = s_{1,1}(n)$ only if $n = 0$.
\end{lemma}

\begin{proof}
    We first observe that
    $$f_{1,1}(z) = \frac{1}{(z-1)^2} = \sum_{n \geq 0} (n+1)z^n,$$
    so for all $n \in \mathbb N$, $$s_{1,1}(n) = n + 1.$$

    In order to prove the inequality, we observe that $s_{1,1}(n)$ is the number of pairs $(X, Y) \in \mathbb N \times \mathbb N$ such that:
    $$X + Y = n,$$
    while $s_{a,b}(n)$ is the number of pairs $(x,y) \in \mathbb N \times \mathbb N$ such that:
    $$ax + by = n.$$
    But $(x,y) \leadsto (X, Y) = (ax, by)$ is a one to one mapping from the set of solutions of the second equation to the set of solutions of the first one. If $a > 1$ this mapping cannot be onto, as its images always have multiples of $a$ on the first component, while for every $0 \leq k \leq n$, the pair $(X, Y) = (k, n-k)$ satisfies the first equation.
\end{proof}

\section{Main results}

\begin{theorem}\label{TheoDivMod}
    Let $a, b \geq 1$ be two natural numbers. Then:
    $$\gcd(a,b) = \left ( \floor{\frac{5^{ ab(ab + a + b)}}{(5^{a^2 b} - 1)(5^{ab^2}-1)}} \bmod 5^{ab} \right ) - 1.$$
\end{theorem}

\begin{proof}
    We apply Lemma \ref{LemmaTermExtraction}. The sequences $s_{a,b}(n)$ have a uniform arithmetic-term representation given by:
    $$\left \lfloor c^{n^2} f_{a,b}(c^{-n}) \right \rfloor \bmod c^n = \left \lfloor \frac{c^{n^2 + an + bn}}{(c^{an} -1)(c^{bn} - 1)} \right \rfloor \bmod c^n,$$
    where a priori the constant $c$ depends on the pair $(a,b)$. We will show that one can find constants  $c \in \mathbb N$ that successfully represent all sequences $s_{a,b}(n)$ for $n \geq m$, where the rank $m$ depends only on $c$.

    Let $R_{a,b}$ be the convergence radius at $z = 0$ for the meromorphic function $f_{a,b}(z)$. The set of poles ${\cal P}(f_{a,b})$ consists of roots of $1$ only:
    $${\cal P}(f_{a,b}) = {\cal U}(a) \cup {\cal U}(b),$$
    where ${\cal U}(m) = \{z \in \mathbb C\,|\, z^m = 1 \}$. As all these poles $u \in \mathbb C$ have $|u| = 1$, we deduce that uniformly $R_{a,b} = 1$ for all $a, b \geq 1$. Hence the condition $c^{-1} < R_{a,b}$ is uniformly satisfied for all $c \in \mathbb N$, $c \geq 2$.

    The condition $s_{a,b} (n) < c^{n-2}$ for all $n \geq m$ and some $m \in \mathbb N$ can be uniformly satisfied for all $a, b \geq 1$. Indeed, by Lemma \ref{LemmaUniformBound}, for $n > 0$ one has:
    $$ s_{a,b}(n) < n+1.$$
    So we observe that already for $c = 2$,
    $$s_{a,b}(n) < n + 1 < 2^{n-2}$$
    for all $n \geq 5$. This allows us to uniformly take $c = 2$ and to be sure that our representation works for $n \geq 5$ for all $a, b \geq 1$. We still have to make the substitution $n = ab$ and get  that:
    $$\gcd(a, b) + 1 = \left \lfloor \frac{c^{ab(ab + a + b)}}{(c^{a^2b} -1)(c^{ab^2} - 1)} \right \rfloor \bmod c^{ab}$$
    for $ab \geq 5$ and $c = 2$. However, as we want a formula which is true for all $a, b \geq 1$, we have to look for bigger values of $c$. The condition that one gets a good result for $a = b = 1$ reads:
    $$2 = \left \lfloor \frac{c^3}{(c -1)^2} \right \rfloor \bmod c.$$
    For $c = 2, 3, 4$ the right-hand side evaluates to $0, 0, 3$ respectively. The first value of $c$ to satisfy this condition is $c = 5$. For $c = 5$, the inequality
    $$n + 1 < 5^{n-2}$$
    is fulfilled by all $n \geq 3$. By direct computation, we observe that the $\gcd$-formula works also for $a, b \geq 1$ with $ab < 3$.
\end{proof}

\begin{corollary}\label{CorModModRep}
    Let $a, b \geq 1$ be two natural numbers. Then:
    $$\gcd(a,b) = \left ( \left ( \left (-5^{ ab(ab + a + b)} \right ) \bmod \left ((5^{a^2 b} - 1)(5^{ab^2}-1) \right ) \right ) \bmod 5^{ab} \right ) - 2.$$
\end{corollary}

\begin{proof} We introduce the following notations: $A:= 5^{ab(ab+a+b)}$, $B : = (5^{a^2 b} - 1)(5^{ab^2}-1)$ and $C:= 5^{ab}$. In order to apply Lemma \ref{LemmaModIdentity}, we check its conditions. The conditions $A, B> 0$, $B \neq 0$,  $C \geq 2$, $C \, |\,A$, $B \nmid A$, $B \bmod C = 1$ are evidently fulfilled. It remains to show that:
$$\lfloor A/B \rfloor \bmod C \neq C-1.$$
But we know that $\lfloor A/B \rfloor \bmod C = \gcd(a,b) + 1$ and $C = 5^{ab}$, so we have to show that $\gcd(a, b) + 1 \neq 5^{ab} - 1$. As $\gcd(a,b) + 1 = s_{a,b}(ab) < ab +1$ by Lemma \ref{LemmaUniformBound}, it will be enough to show that:
$$ab + 1 < 5^{ab} - 1$$
for all $a, b \geq 1$. Denote $ab = k \geq 1$. We observe that the inequality:
$$k + 2 < 5^k$$
holds for all $k \geq 1$.
\end{proof}

By proving inequalities and verifying the finitely many remaining cases, as in the proofs of Theorem \ref{TheoDivMod} and of Corollary \ref{CorModModRep}, one shows:

\begin{corollary}\label{CorBase2}
    For all natural numbers $a, b \geq 1$ with $(a,b) \neq (1,1)$ the following identities hold true:
    $$\gcd(a,b) = \left ( \floor{\frac{2^{ ab(ab + a + b)}}{(2^{a^2 b} - 1)(2^{ab^2}-1)}} \bmod 2^{ab} \right ) - 1,$$
    $$\gcd(a,b) = \left ( \left ( \left (-2^{ ab(ab + a + b)} \right ) \bmod \left ((2^{a^2 b} - 1)(2^{ab^2}-1) \right ) \right ) \bmod 2^{ab} \right ) - 2.$$
    If one replaces the exponentiation base $c=2$ with the bases $c=3$ or $c=4$, the identities hold true for all $a, b \geq 1$, with the unique exception $(a,b) = (1,1)$ again.
\end{corollary}

\begin{remark}\rm The mod - mod representations given in Corollary \ref{CorModModRep} and in  Corollary \ref{CorBase2} are not arithmetic terms as defined in \cite{mazzanti2002plainbases} and in  \cite{marchenkov2007superposition}, because the subterm $-c^{ab(ab+a+b)}$ is essentially negative. However, the mod - mod representations have computational advantages over the div - mod representations, because the subterm:
$$\left ( \left  (-c^{ ab(ab + a + b)} \right ) \bmod \left ( (c^{a^2 b} - 1)(c^{ab^2}-1) \right ) \right )$$
can be computed much faster than the analogous
$$ \floor{\frac{c^{ ab(ab + a + b)}}{(c^{a^2 b} - 1)(c^{ab^2}-1)}}$$
by using the Algorithm for Fast Exponentiation in modular arithmetic. \qed
\end{remark}


\end{document}